\documentclass[12pt]{article}
\usepackage{amsfonts,amsmath}
\usepackage{amssymb,latexsym}
\usepackage[dvips]{epsfig}
\usepackage{psfrag}
\usepackage{amscd}
\usepackage{pstricks,pst-node}
\usepackage{rotating}
\usepackage{mfpic}

\title{An invariant of tangle cobordisms via subquotients of arc rings}

\author{Yanfeng Chen and Mikhail Khovanov}
\date{October 3, 2006}

\newtheorem{prop}{Proposition}
\newtheorem{theorem}{Theorem}
\newtheorem{lemma}{Lemma}

\newcommand{\Z}{\mathbb{Z}}
\newcommand{\C}{\mathbb{C}}
\newcommand{\R}{\mathbb{R}}
\newcommand{\Uqsl}{U_q(sl(2))}
\newcommand{\oplusop}[1]{{\mathop{\oplus}\limits_{#1}}}
\newcommand{\oplusoop}[2]{{\mathop{\oplus}\limits_{#1}^{#2}}}
\newcommand{\mo}{\mathbf{1}}
\newcommand{\cA}{{\mathcal{A}}}
\newcommand{\cF}{{\mathcal{F}}}
\newcommand{\cK}{{\mathcal{K}}}
\newcommand{\Inv}{\textrm{Inv}}
\newcommand{\define}{\stackrel{\mbox{\scriptsize{def}}}{=}}
\def\doublemaprights#1#2#3#4{\raise3pt\hbox{$\mathop{\,\,\hbox to
     #1pt{\rightarrowfill}\kern-30pt\lower3.95pt\hbox to
     #2pt{\leftarrowfill}\,\,}\limits_{#3}^{#4}$}}

\begin{document}

\maketitle
\baselineskip 14pt

\tableofcontents

\section{Introduction}

The extension of the Jones polynomial of links \cite{Jones} to
tangles is governed, from the algebraic viewpoint, by the quantum
group $\Uqsl$ and its representation theory. In one possible
extension, to $n$ points on the plane there is assigned
$V^{\otimes n},$ the $n$-th tensor power of the fundamental
two-dimensional representation of $\Uqsl,$ and to an
$(m,n)$-tangle $T$ a homomorphism $f(T)$ of representations
$V^{\otimes n} \longrightarrow V^{\otimes m}.$ An $(m,n)$-tangle
is a tangle in $\R^2 \times [0,1]$ with $m$ top and $n$ bottom
endpoints. Alternatively, it's possible to restrict to tangles
with even number of top and bottom endpoints, assign to $2n$
points on the plane the space $\Inv(V^{\otimes 2n})$ of
$\Uqsl$-invariants in $V^{\otimes 2n},$ and to a $(2m, 2n)$-tangle
$T$ the map
$$ f_{inv}: \mathrm{Inv}(V^{\otimes 2n}) \longrightarrow
   \mathrm{Inv}(V^{\otimes 2m}),$$
the restriction of $f$ to the subspace of invariants. When the
tangle is a link, $f(T)=f_{inv}(T)$ is the endomorphism of a one-dimensional
vector space, given by multiplication by the Jones polynomial of $T.$

A categorification of the Jones polynomial \cite{Kh1} was
extended  to tangles and tangle cobordisms in \cite{Kh2}, \cite{Kh3}.
The space $\Inv(V^{\otimes 2n}),$
interpreted as a free $\Z[q,q^{-1}]$-module of rank equal to the $n$-th
Catalan number, becomes the Grothendieck group of the triangulated
category $\mathcal{K}_n,$ which is the category of bounded complexes
of finitely-generated graded modules over a certain graded ring $H^n,$
up to chain homotopies of complexes. The invariant $\cF(T)$ of a tangle $T$
is an appropriate exact functor $\mathcal{K}_n \longrightarrow \mathcal{K}_m,$
which on the Grothendieck group gives the map $f_{inv}(T).$ The invariant
of a tangle cobordism is a natural transformation of functors.

Bar-Natan \cite{BN} suggested a more general and geometric
categorification of the Jones polynomial and its extension to
tangles and their cobordisms (also see \cite{N}). From the
algebraic viewpoint, he considers the universal deformation of the
original construction, with homology defined over a bigger ground
ring, and tangle invariants taking values in a category similar
but richer than $\mathcal{K}_n.$ The Grothendieck group of his
category is, again, naturally isomorphic to a
$\Z[q,q^{-1}]$-lattice in $\mathrm{Inv}(V^{\otimes n}).$

In each of these two examples, generalization of the
link homology to tangles utilizes a categorification of the space of
invariants  $\mathrm{Inv}(V^{\otimes 2n}).$
The categorification \cite{KR} of the quantum $sl(m)$ invariant of links and tangles,
when specialized to $m=2,$ uses categories of complexes of matrix factorizations with
potentials $\sum \pm x_i^3,$ which contain proper subcategories equivalent to
$\mathcal{K}_n$ (over $\mathbb{Q}$). The Grothendieck groups
of these categories of matrix factorizations have not been computed.

A categorification of the entire tensor product $V^{\otimes n}$ was investigated
in \cite{BFK}. One first forms  a suitable direct sum
$$ \mathcal{O}^n = \oplusoop{k=0}{n} \mathcal{O}^{n-k,k} $$
of parabolic blocks of the highest weight category for $sl(n).$
The category $\mathcal{O}^{n-k,k}$ is a full subcategory of a
regular block of $\mathcal{O}$ for $sl(n)$ consisting of modules
which are locally-finite as $sl(n-k)\times sl(k)$-modules. The
Grothendieck group of $\mathcal{O}^{k,n-k}$ is free abelian of
rank $\left( \begin{array}{c} n \\ k \end{array}
 \right)$ and, after tensoring with $\C$ over $\Z,$ can
be naturally identified with the weight $n-2k$ subspace in
$V^{\otimes n}$
  $$ K(\mathcal{O}^{n-k,k}) \otimes_{\Z} \C \cong V^{\otimes n} (n-2k).$$
 The action of the Temperley-Lieb algebra on $V^{\otimes n}$ lifts to
exact endofunctors (called translation across the wall) in $\mathcal{O}^n.$

The next and major development in this direction was due to
Stroppel \cite{St1}, \cite{St2}, who considered a graded version of $\mathcal{O}^n$
and of the translation functors, and showed that they produce an invariant of
tangles and tangle cobordisms. In this extension the invariant of an $(m,n)$-tangle
$T$ is a functor $D^b(\mathcal{O}^n) \longrightarrow
D^b(\mathcal{O}^m)$ between the derived categories.

The  $U_q(sl(2))$-invariants $\mathrm{Inv}(V^{\otimes 2n})$ form a subspace of
the weight zero space $V^{\otimes 2n}(0)$ of $V^{\otimes 2n}.$
One would expect that the inclusion
$$\mathrm{Inv}(V^{\otimes 2n}) \subset V^{\otimes 2n}(0)$$
can be lifted to the level of categories, to some relation between
$\mathcal{K}_n$ and $D^b(\mathcal{O}^{n,n}).$ That this is indeed
the case was conjectured by Stroppel \cite{St2}  and recently
proved by her in \cite{St3}. Namely, the category
$\mathcal{O}^{n-k,k}$ is equivalent to the category of
finite-dimensional modules over certain finite-dimensional basic
$\C$-algebra $A_{n-k,k},$ described by Braden \cite{Braden} via
generators and relations. Stroppel showed that the ring $H^n$
which controls the categorification $\mathcal{K}_n$ of the
invariant space is isomorphic to a subring of $ A_{n,n}.$ More
precisely,
$$H^n\otimes_{\Z} \C\cong e A_{n,n}e, $$
where $e$ is an idempotent such that $A_{n,n}e$ is the largest direct
summand of $A_{n,n}$ which is both a projective and an injective $A_{n,n}$-module,
see \cite{St3}.

Our paper was motivated by the problem of categorifying
$V^{\otimes n}$ and the linear maps $f(T)$ directly, in a
down-to-earth way, avoiding the highly sophisticated machinery of
highest weight categories and their graded versions. We define a
collection of graded rings $A^{n-k,k},$ consider the categories of
finitely-generated graded $A^{n-k,k}$-modules, and identify their
Grothendieck groups with $\Z[q,q^{-1}]$-lattices in the weight
spaces of $V^{\otimes n}.$ Next, we form product rings
$$ A^n \define \prod_{k=0}^n A^{n-k,k}, $$
to an $(m,n)$-tangle $T$ assign a complex $\mathcal{F}(T)$
of graded $(A^m, A^n)$-bimodules, and to a tangle cobordism--a homomorphism
of complexes.

During our work on this project, Stroppel's paper \cite{St3} came out, where
she defines a ring $\mathcal{K}^n$ isomorphic to $A^{n,n}\otimes \C$
and shows that $\mathcal{K}^n$ is isomorphic to the Braden algebra
$A_{n,n}$ \cite{Braden}. Furthermore,
Stroppel announced the theorem that the inclusion of subrings
$H^n\otimes_{\Z}\C\subset A_{n,n}$ extends to bimodules and bimodule
homomorphisms in the two theories, allowing her to directly relate
tangle and tangle cobordism invariants of \cite{Kh2}, \cite{Kh3} with
those of \cite{St1, St2}.

Our constructions and results have a nonempty intersection
with Stroppel \cite{St3}, and, as we expect, will be easily surpassed by
her announced work. We decided to publish this paper, nevertheless,
since our work, which was done independently, will be a basis for the
forthcoming paper \cite{YC}.

\vspace{0.1in}

{\bf Acknowledgements:} M.K. was partially supported by the NSF grant DMS-0407784.
Y.C. would like to thank his advisor Robion Kirby for encouragement and
financial support.

\section{Arc ring $H^n$}
We first recall the definition of $H^n$ from
\cite{Kh2}. Let $\cA$ be a free graded
abelian group of rank $2$ spanned by $\mo$ in degree $-1$ and $X$
in degree $1$. Define the unit map $\iota: \mathbb{Z} \rightarrow
\cA$ and the trace map $\epsilon: \cA \rightarrow \mathbb{Z}$ by
 \begin{equation*}
 \iota(1)=\mo,\hspace{0.1in} \varepsilon(\mo)=0, \hspace{0.1in}\varepsilon(X)=1.
 \end{equation*}
Define multiplication $m: \cA \otimes \cA \rightarrow \cA$ by
 \begin{equation}
 \mo^2=\mo, \hspace{0.1in} \mo X= X\mo = X, \hspace{0.1in} X^2=0, \label{m}
 \end{equation}
and comultiplication $\Delta$ by
 \begin{equation}
 \Delta: \cA \rightarrow \cA^{\otimes 2}, \hspace{0.1in}\Delta(\mo) = \mo \otimes X + X\otimes \mo , \hspace{0.1in} \Delta(X) = X\otimes
 X.\label{Delta}
 \end{equation}

\begin{figure}[ht!]
\begin{center}
\psfrag{s01}{$S^1_0$}
\psfrag{s10}{$S^0_1$}\psfrag{s21}{$S^1_2$}\psfrag{s12}{$S^2_1$}
\epsfig{figure=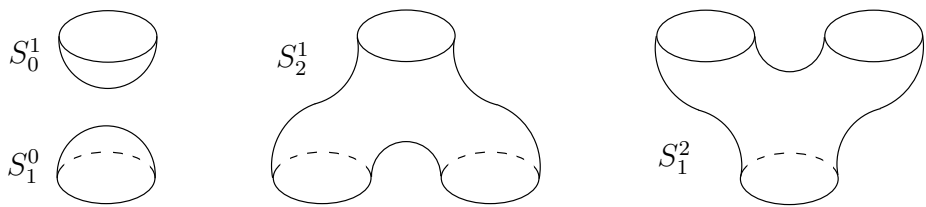} \caption{Elementary
cobordisms.} \label{Elementry Cobs.figure}
\end{center}
\end{figure}

Assign to $\cA$ a $2$-dimensional TQFT $\cF$ which associates
$\cA^{\otimes k}$ to a disjoint union of $k$ circles. To the
elementary cobordisms $S_0^1$, $S_1^0$, $S_2^1$ and $S_1^2$,
depicted in figure~\ref{Elementry Cobs.figure}, $\cF$ associates
maps $\iota$, $\varepsilon$, $m$ and $\Delta$ respectively. The
map $\mathcal{F}(S)$ between tensor powers of $\cA$ induced by a
cobordism $S$ is homogeneous of degree minus the Euler
characteristic of $S$
 \begin{equation}deg(\mathcal{F}(S))=-\chi(S).\label{DegreeOfCobordism}
 \end{equation}

\begin{figure}[ht!]
\begin{center}
\epsfig{figure=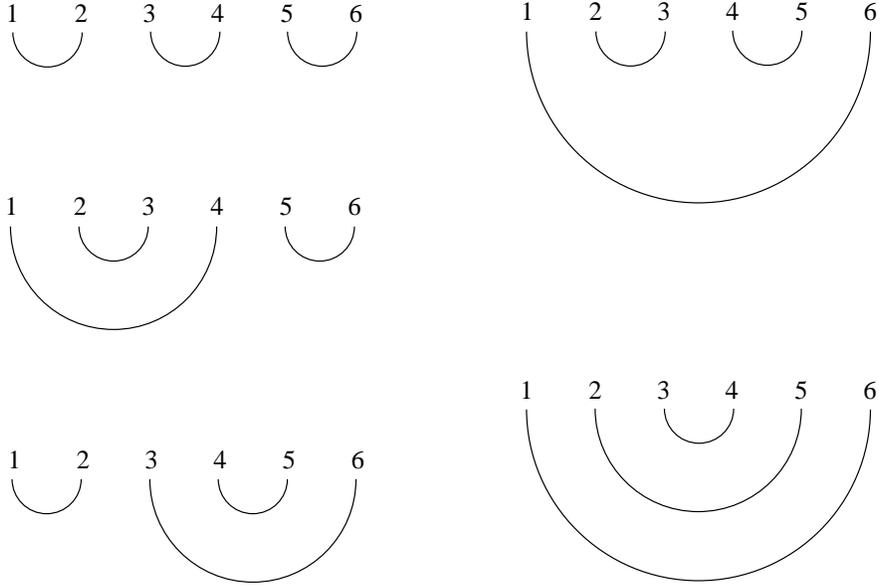} \caption{Crossingless matchings of $6$
points.} \label{B3.figure}
\end{center}
\end{figure}

\begin{figure}[ht!]
\begin{center}
\psfrag{b}{\small$b$}\psfrag{a}{\small$a$}\psfrag{W(b)}{\small$W(b)$}\psfrag{W(b)a}{\small$W(b)a$}

\epsfig{figure=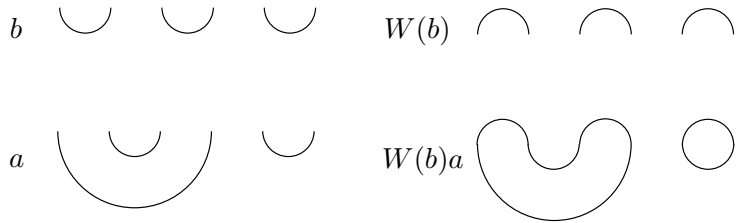}  \caption{Gluing in $B^3$.}
\label{Wba.figure}
\end{center}
\end{figure}

Let $B^n$ be the set of crossingless matchings of $2n$ points.
figure~\ref{B3.figure} shows the set $B^3$.
For $a$, $b\in B^n$ denote by $W(b)$ the reflection of $b$ about
the horizontal axis, and by $W(b)a$ the closed $1$-manifold
obtained by closing $W(b)$ and $a$ along their boundaries, see
figure~\ref{Wba.figure}.

The graded abelian group $\cF(W(b)a)$ is isomorphic to
$\cA^{\otimes I}$, where $I$ is the set of circles in $W(b)a$. The
symbol $\{n\}$ denotes shifting the grading up by $n$.
For $a$, $b\in B^n$ let
 \begin{equation*}
 {_b(H^n)_a} \define \cF(W(b)a)\{n\},
 \end{equation*}
and define $H^n$ as the direct sum
 \begin{equation*}
 H^n\define \oplusop{a,b\in B^n} \hspace{0.05in} {_b(H^n)_a}
 \end{equation*}

\begin{figure}[ht!]
\begin{center}
\psfrag{a}{\tiny$a$}\psfrag{c}{\tiny$c$}\psfrag{W(b)}{\tiny$W(b)$}\psfrag{W(a)}{\tiny$W(a)$}

\epsfig{figure=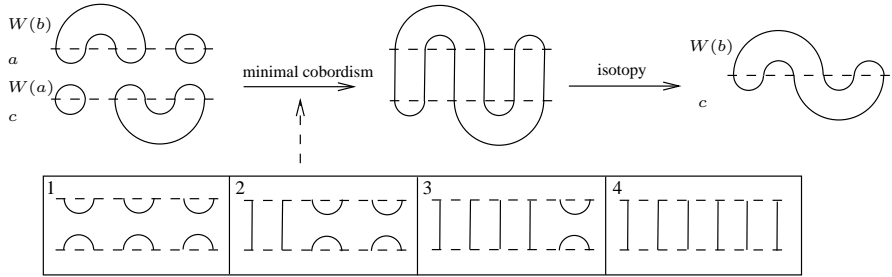} \caption{Multiplication in $H^n$.}
\label{contraction.figure}
\end{center}
\end{figure}

 Multiplication maps in
$H^n$ is defined as follows. We set $xy=0$ if $x\in {_b(H^n)_a}$,
$y\in {_c(H^n)_d}$ and $c\neq a$. Multiplication maps
 \begin{equation*}
 {_b(H^n)_a} \otimes {_a(H^n)_c} \rightarrow {_b(H^n)_c}
 \end{equation*}
are given by homomorphisms of abelian groups
 \begin{equation*}
 \cF(W(b)a) \otimes \cF(W(a)c) \rightarrow \cF(W(b)c)
 \end{equation*}
which are induced by ``minimal'' cobordisms from $W(b)aW(a)c$ to
$W(b)c$, see figure~\ref{contraction.figure}.

The element $1_a\define \mo^{\otimes n}\in \cA^{\otimes n} \cong
{_a(H^n)_a}$ is an idempotent in $H^n.$ The sum $\sum_a 1_a$ is
the unit element of $H^n$. See \cite{Kh2} for details.

\section{Subquotients of $H^{n}$}
\begin{figure}[ht!]
\begin{center}
\epsfig{figure=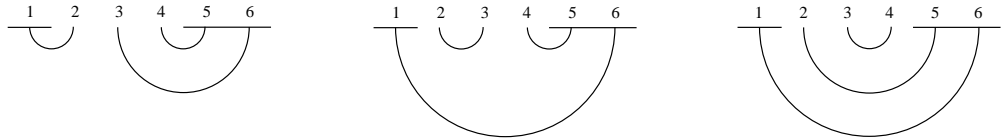}  \caption{The $3$ elements in $B^{1,2}$.}
\label{b21.figure}
\end{center}
\end{figure}

For each $n\geq 0$ and $0\leq k \leq n$, define $B^{n-k,k}$ to be
the subset of $B^n$ consisting of diagrams with no matchings among
the first $n-k$ points and among the last $k$ points.
Figure~\ref{b21.figure} shows $B^{1,2}$ (compare with
figure~\ref{B3.figure}).
We put two ``platforms'', one on the first $n-k$ points and one on
the last $k$ points, to indicate that these endpoints are special.
Define $\widetilde{A}^{n-k,k}$ by
 \begin{equation}
   \widetilde{A}^{n-k,k} \define \oplusop{a,b\in B^{n-k,k}}
   \hspace{0.05in} \cF(W(b)a)\{n\}. \label{DefTildeA.equation}
 \end{equation}
$\widetilde{A}^{n-k,k}$ sits inside $H^n$ as a graded subring
which inherits its multiplication from $H^n$ (the inclusion takes
$1\in \widetilde{A}^{n-k,k}$
to an idempotent of $H^n$).\\

\begin{figure}[ht!]
\begin{center}
\psfrag{1}{\tiny I} \psfrag{2}{\tiny II} \psfrag{3}{\tiny III}
\epsfig{figure=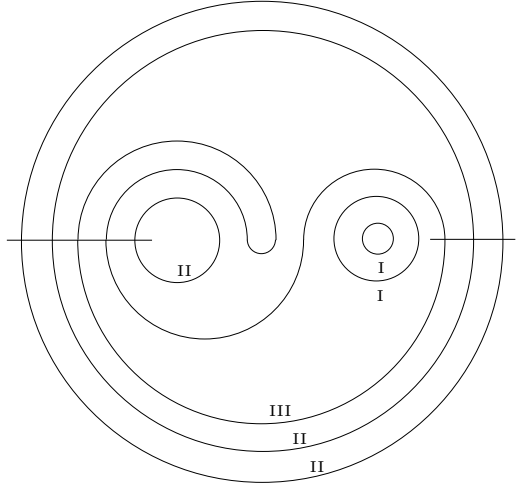} \caption{$3$ types of circles in
$\cF(W(b)a)$.} \label{3types.figure}
\end{center}
\end{figure}

For $a,b\in B^{n-k,k}$ the circles of $W(b)a$ fall into $3$
different types (see figure~\ref{3types.figure}):
\begin{itemize}
\item Type I: Circles that are disjoint from platforms.

\item Type II: Circles that intersect at least one platform and
intersect each platform at most once.

\item Type III: Circles that intersect one of the platforms at
least twice.
\end{itemize}

\begin{figure}[ht!]
\begin{center}
\epsfig{figure=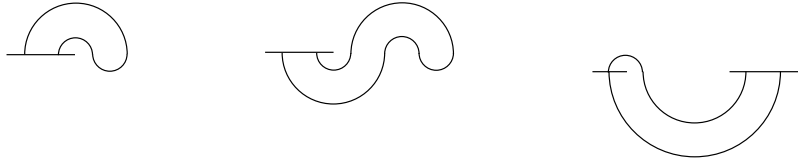} \caption{Examples of portions of
type III circles.} \label{ZeroElement.figure}
\end{center}
\end{figure}

We call an intersection point between a circle and a platform a
``mark''. Next, we introduce an ideal $I^{n-k,k}\subset
\widetilde{A}^{n-k,k}$. If $W(b)a$ contains at least one type III
circle (see figure~\ref{ZeroElement.figure}), set
${_b(I^{n-k,k})_a}=\cF(W(b)a).$  If $W(b)a$ contains only circles
of type I and II, we write $\cF(W(b)a)=\cA^{\otimes i}\otimes
\cA^{\otimes j},$ where type II circles correspond to the first
$i$ tensor factors, and define ${_b(I^{n-k,k})_a}$ as the span of
 \begin{equation*}
y_1 \otimes \cdots \otimes y_{t-1} \otimes X \otimes y_{t+1}
\otimes \cdots \otimes y_{i+j} \in   \cA^{\otimes i}\otimes
\cA^{\otimes j} \cong \cF(W(b)a),
 \end{equation*}
where $1\leq t\leq i$ and $y_s\in \{\mo, X\}$. By taking the
direct sum over all $a,b\in B^{n-k,k}$ we get a subgroup of
$\widetilde{A}^{n-k,k}$

 \begin{equation*}
   I^{n-k,k} \define \oplusop{a,b\in B^{n-k,k}}\hspace{0.05in}
   {_b(I^{n-k,k})_a}.
 \end{equation*}

\begin{figure}[ht!]
\begin{center}
\psfrag{1}{\small $\mo$}\psfrag{x}{\small $X$}
\epsfig{figure=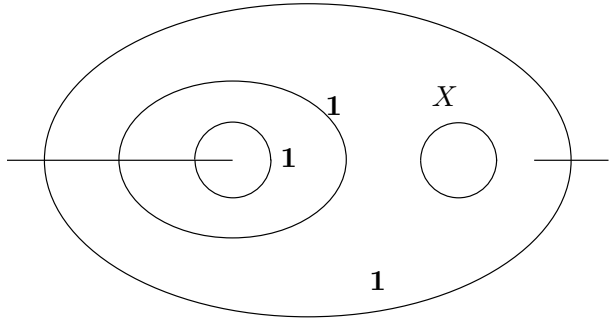} \caption{Visualization of a
tensor product element in $\cF(W(a)b)$.}
\label{VisualizeElements.figure}
\end{center}
\end{figure}

\begin{lemma} \label{lemma-ideal}
$I^{n-k,k}$ is a two-sided graded ideal of the ring
$\widetilde{A}^{n-k,k}$.
\end{lemma}
\emph{Proof:} To prove it's a left ideal, it suffices to show that
$uv\in {_c(I^{n-k,k})_b}$ whenever $u\in \cF(W(c)a)$ and $v\in
{_a(I^{n-k,k})_b}$. Without loss of generality, we can assume that
both $u$ and $v$ are tensor products
 \begin{equation*}
u=u_1 \otimes \cdots \otimes u_{s} \in   \cA^{\otimes s} \cong
\cF(W(c)a),
 \end{equation*}
 and
  \begin{equation*}
v=v_1 \otimes \cdots \otimes v_{t} \in   \cA^{\otimes t} \cong
\cF(W(a)b),
 \end{equation*}
where $u_i$, $v_j\in \{\mo, X\}$. We can visualize $u$ and $v$ as
sets of circles with labels $\mo$ or $X$, see
figure~\ref{VisualizeElements.figure}.

Case 1: ${_a(I^{n-k,k})_b}\neq \cF(W(a)b)$. In this case, $v$
contains a marked circle $C$ with label $X$. Pick a mark $p$ on
$C$. Denote by $M\in \{\mo, X\}$ the label of the circle
containing $p$ in $W(c)b$. It follows from equations (\ref{m}) and
(\ref{Delta}) that after multiplying $v$ by $u$, either $M$ will
remain $X$ or $uv=0$. In either case $uv$ will belong to
${_c(I^{n-k,k})_b}$. See figure~\ref{ideal.figure} for a similar
example.

Case 2: ${_a(I^{n-k,k})_b}= \cF(W(a)b)$. In this case $v$ contains
a circle connecting two points $p_1$ and $p_2$ on the same
platform. If $C$ is labelled by $X$, it follows from the previous
case that either the labels of circles containing $p_1$ and $p_2$
will remain $X$ or $uv=0$. So $uv$ will belong to
${_c(I^{n-k,k})_b}$. Now assume that the label on $C$ is $\mo$.
If, during the process of multiplying $u$ and $v$, a splitting of
$C$ takes place it follows from equation (\ref{Delta}) that either
the circle containing $p_1$ or the circle containing $p_2$ will
have label $X$. Otherwise, a sequence of merging with $C$ will
keep $p_1$ and $p_2$ connected by a single arc. In this case
${_c(I^{n-k,k})_b}=\cF(W(c)b)$, and therefore contains $uv$.

Similar arguments show that $I^{n-k,k}$ is a right ideal, and the
lemma follows. $\square$ \vspace{0.1in}

 The ring $A^{n-k,k}$ is defined as the
quotient of $\widetilde{A}^{n-k,k}$ by the ideal $I^{n-k,k}$
 \begin{equation}
   A^{n-k,k} \define \widetilde{A}^{n-k,k} / I^{n-k,k}. \label{DefA.equation}
 \end{equation}
 $A^{n-k,k}$ naturally decomposes into a direct sum of graded
 abelian groups
 \begin{equation*}
   A^{n-k,k} =  \oplusop{a,b\in B^{n-k,k}} \hspace{0.05in}
   _a(A^{n-k,k})_b,
 \end{equation*}
where $_a(A^{n-k,k})_b = \cF(W(a)b) /
{_a(I^{n-k,k})_b}\{n\}$. The abelian group ${_a(A^{n-k,k})_b}=0$ if
and only if $W(b)a$ contains a type III circle. Otherwise,
${_a(A^{n-k,k})_b}$ is a free abelian group of rank $2^{c_1}$
where $c_1$ is the number of type I circles in $W(b)a$. Assuming
that $\cF(W(a)b)\cong {\cA}^{\otimes m}$ in which type II circles
correspond to the first $i$ tensor factors, ${_a(A^{n-k,k})_b}$
has a basis of the form
 \begin{equation*}
\mo \otimes \cdots \otimes \mo \otimes a_{i+1} \otimes \cdots
\otimes a_{m}
 \end{equation*}
 where $a_s\in \{\mo, X\}$ for all $i+1\leq s\leq m$.\\

The element $1_a\define \mo^{\otimes n}\in {_a(A^{n-k,k})_a}$ is a
minimal idempotent in $A^{n-k,k}.$ The sum $1 \define \sum_{a\in
B^{n-k,k}} 1_a$ is the unit element of $A^{n-k,k}$.

The relations among the three rings $H^n$,
$\widetilde{A}^{n-k,k}$, and $A^{n-k,k}$ is described in the
following diagram:

\newpage

\begin{figure}[ht!]
\begin{center}
\psfrag{h}{$H^n$}\psfrag{a1}{$\widetilde{A}^{n-k,k}$}\psfrag{a2}{$A^{n-k,k}$}
\psfrag{s}{\small Inclusion of subrings}\psfrag{q}{\small Quotient
by $I^{n-k,k}$}
 \epsfig{figure=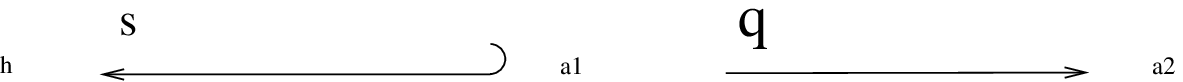}
\end{center}
\end{figure}

We now write down explicitly the rings $A^{n-k,k}$ in simplest
cases.
\begin{itemize}
\item The ring $A^{0,n}$ is isomorphic to $\mathbb{Z}$, since
$B^{0,n}$ contains only one diagram, and the functor $\cF$ applied
to its closure produces $\mathbb{Z}$ (see
figure~\ref{b0n.figure}).

\begin{figure}[ht!]
\begin{center}
\psfrag{1}{\tiny $1$} \psfrag{n-1}{\tiny $n-1$} \psfrag{n}{\tiny
$n$} \psfrag{n+1}{\tiny $n+1$} \psfrag{n+2}{\tiny $n+2$}
\psfrag{2n}{\tiny $2n$}\psfrag{cd}{\tiny
$\cdots$}\psfrag{a}{$a$}\psfrag{W}{$W(a)a$}
 \epsfig{figure=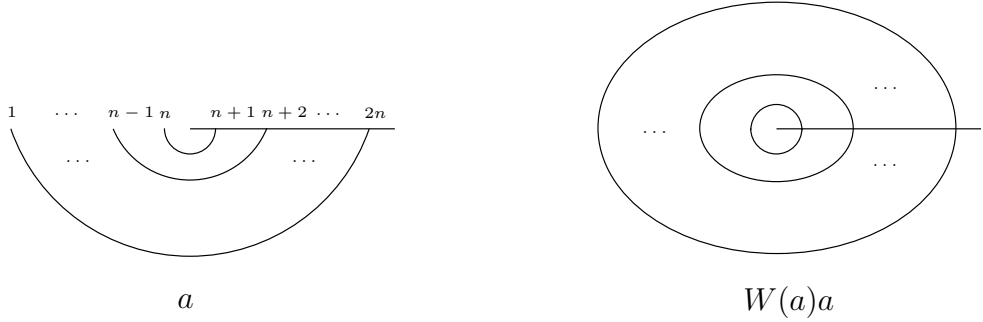}  \caption{The only element in $B^{0,n}$ and its closure.} \label{b0n.figure}
\end{center}
\end{figure}

\item The vertical reflection induces an isomorphism $A^{n-k,k}\cong A^{k,n-k}.$ Namely,
 reflecting a diagram in $B^{n-k,k}$ about a vertical axis produces a diagram
in $B^{k,n-k},$ leading to an isomorphism of sets $B^{n-k,k}\cong
B^{k,n-k}$. This isomorphism induces an isomorphism of rings
$\widetilde{A}^{n-k,k}\cong \widetilde{A}^{k,n-k}$ and of the
quotient rings $A^{n-k,k}\cong A^{k,n-k}$. In particular,
$$A^{n,0}\cong A^{0,n}\cong \mathbb{Z}.$$

\item The set $B^{1,1}$ contains two diagrams which we denote by
$a$ and $b$ respectively (see figure~\ref{b11.figure}).

\begin{figure}[ht!]
\begin{center}
\psfrag{a}{$a$} \psfrag{b}{$b$} \epsfig{figure=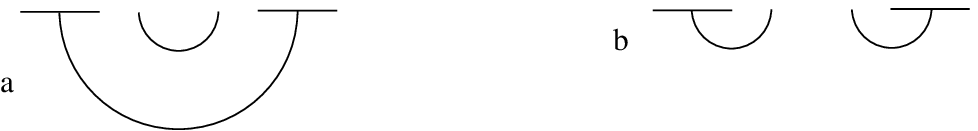}
\caption{$B^{1,1}$.} \label{b11.figure}
\end{center}
\end{figure}

From figure~\ref{a11.figure} we can see that
\[
\begin{array}{rclcrcl}
_a(A^{1,1})_a & = & \cA\{1\}, & \hspace{0.2in} &
_b(A^{1,1})_a & = & \mathbb{Z}\{1\}, \\
_a(A^{1,1})_b & = & \mathbb{Z}\{1\},  & \hspace{0.2in} &
_b(A^{1,1})_b & = & \mathbb{Z},
\end{array}
\]
where $\{1\}$ denotes shifting the grading up by $1$. The grading
shifts in above formulas follow from the definition of
$A^{n-k,k}$. For example,
${_b(A^{1,1})_b}=\cF(W(b)b)/{_b(I^{1,1})_b}\{2\}$ generated by
$\mo\otimes\mo\{2\}$ which sits in degree $0.$ The ring
$A^{1,1}$ has a simple quiver description, as the path ring of the graph
\begin{equation*}
\label{sl2}
 \stackrel{a}{\circ} \doublemaprights{30}{30}{\beta}{\alpha} \stackrel{b}{\circ}
\end{equation*}
with the defining relation $\alpha \beta=0$. Paths $\alpha$ and
$\beta$ correspond to generators of ${_b(A^{1,1})_a}$ and
${_a(A^{1,1})_b}$ respectively.
\begin{figure}[ht!]
\begin{center}
\psfrag{waa}{$W(a)a$}
\psfrag{wab}{$W(a)b$}\psfrag{wba}{$W(b)a$}\psfrag{wbb}{$W(b)b$}
\epsfig{figure=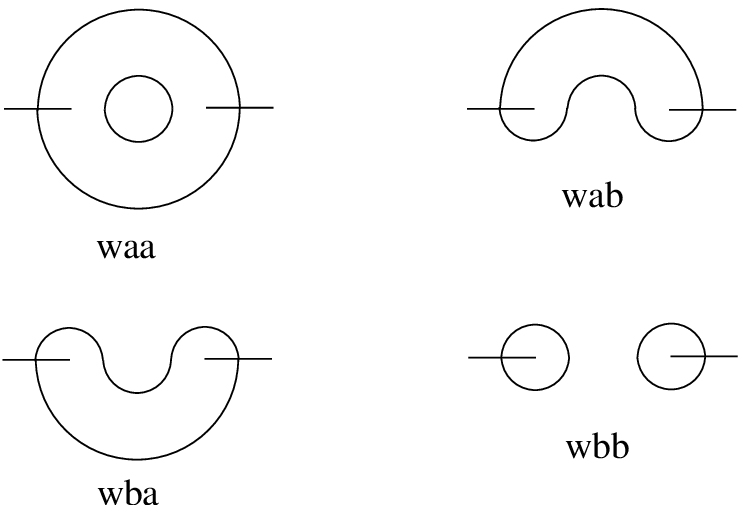}  \caption{$A^{1,1}$.} \label{a11.figure}
\end{center}
\end{figure}
\item The ring $A^{1,n-1}$ has a simple quiver description as
well. Denote by $Q_n$ the quotient of the path ring
\begin{equation*}
\label{sl20}
 \stackrel{1}{\circ} \doublemaprights{30}{30}{}{} \stackrel{2}{\circ}
  \doublemaprights{30}{30}{}{} \stackrel{3}{\circ}
 \doublemaprights{30}{30}{}{} \dots
   \doublemaprights{30}{30}{}{}
\stackrel{n-1}{\circ}
 \doublemaprights{30}{30}{}{} \stackrel{n}{\circ}
\end{equation*}
by the defining relations $(1|2|1)=0$, $(i|i+1|i)= (i|i-1|i)$ for
$1<i<n$, $(i|i+1|i+2)=0$ for $i<n-1$, and $(i|i-1|i-2)=0$ for
$i\geq 2$, where $(i|j|k)$ denotes the path which starts at $i$, goes to
$j$ and then to $k$. See \cite{KS} for a more detailed
discussion of this algebra.

The rings  $Q_n$ and $A^{1,n-1}$ are isomorphic. The isomorphism
takes the idempotent $1_{a_i}$, for $a_i$ shown in
figure~\ref{b13.figure}, to the minimal idempotent $(i)$, which is
the length zero path that starts and ends at $i$. The path $(j|i)$
with $j=i\pm 1$ corresponds to a generator of $\cF(W(a_j)a_i)\cong
\mathbb{Z}$.

Figure~\ref{b13.figure} depicts all $a_i$ for $A^{1,3}$. The
diagram $W(a_3)a_1$ has a type III circle so
${_{a_3}(A^{1,3})_{a_1}}=0$. More generally,
${_{a_i}(A^{1,3})_{a_j}}=0$ iff $|i-j|>1$. See
figure~\ref{a13.figure} for an example of the multiplication in
$A^{1,3},$ and note that $Q_4$ is given by
\begin{equation*}
\label{sl21}
 \stackrel{1}{\circ} \doublemaprights{30}{30}{}{} \stackrel{2}{\circ}
  \doublemaprights{30}{30}{}{} \stackrel{3}{\circ}
 \doublemaprights{30}{30}{}{} \stackrel{4}{\circ}
\end{equation*}
with the relations $(1|2|1)=0$, $(2|1|2)= (2|3|2)$, $(3|2|3)=
(3|4|3)$, and $(1|2|3)=(2|3|4)=(4|3|2)=(3|2|1)=0$.

\begin{figure}[ht!]
\begin{center}
\psfrag{a1}{$a_1$}\psfrag{a2}{$a_2$}
\psfrag{a3}{$a_3$}\psfrag{a4}{$a_4$} \psfrag{1}{\small $(1)$}
\psfrag{2}{\small $(2)$}\psfrag{3}{\small $(3)$}\psfrag{4}{\small
$(4)$} \epsfig{figure=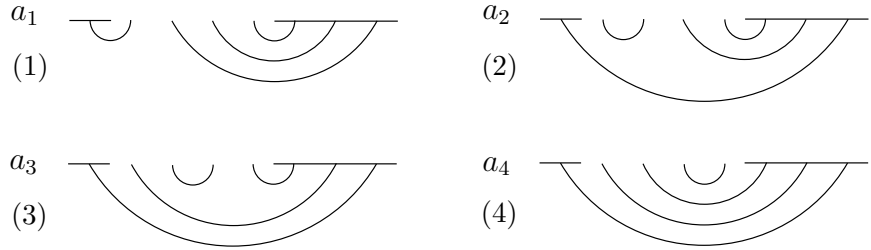} \caption{Elements in
$B^{1,3}$.} \label{b13.figure}
\end{center}
\end{figure}

\begin{figure}[ht!]
\begin{center}
\psfrag{a12}{$W(a_2)a_1$}\psfrag{a23}{$W(a_3)a_2$}\psfrag{a13}{$W(a_3)a_1$}
\psfrag{12}{\small $(1|2)$}\psfrag{23}{\small
$(2|3)$}\psfrag{13}{\small $(1|3)$} \psfrag{m}{\tiny
minimal}\psfrag{c}{\tiny cobordism} \epsfig{figure=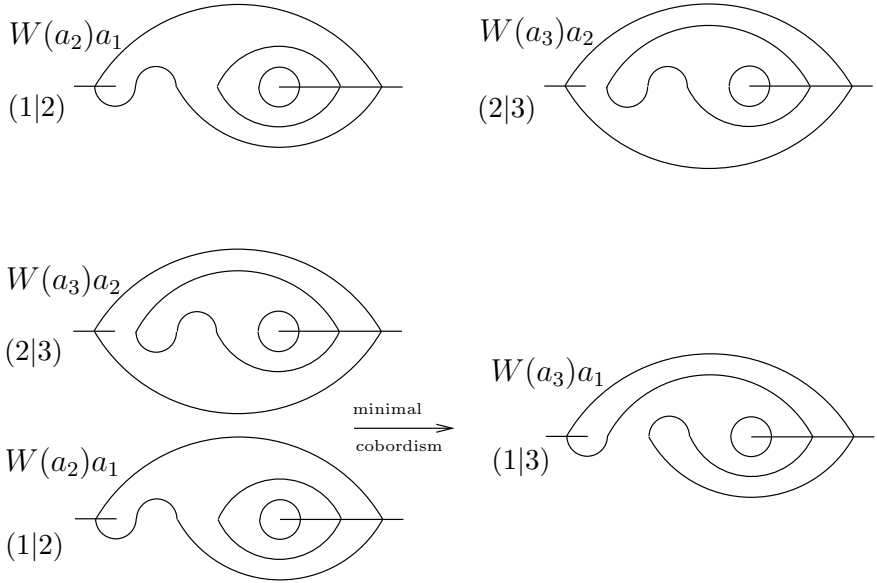}
\caption{Quiver description of $A^{1,3}$.} \label{a13.figure}
\end{center}
\end{figure}

\end{itemize}

We list some basic properties of rings $A^{n-k,k}$:

\begin{itemize}

\item The rings $A^{n-k,k}$ are indecomposable ($0$ and $1$ are
the only central idempotents in $A^{n-k,k}$).

\item $B^{n-k,k}$ is the union of two disjoint subsets $B_1$ and
$B_2$ as follows. Any two elements in $B^{k,n-k}$ can be connected
by a sequence of elementary changes as in
figure~\ref{ElementaryMove.figure}. Pick an element $a$ in
$B^{n-k,k}$ and put it in $B_1$. For any $b\in B^{n-k,k}$, if $a$
and $b$ can be connected in an even number of steps we put $b$ in
$B_1$, otherwise put it in $B_2$. By taking the sum of all minimal
idempotents in each subset we get idempotents
$$e_1=\sum_{a\in B_1}1_a,\ \ \ \ \ \  e_2=1-e_1=\sum_{a\in B_2}1_a,$$
such that all homogeneous elements in $e_i A^{n-k,k} e_i$, for
$i$, $j$ in the same subset, have even degrees, and all
homogeneous elements in $e_i A^{n-k,k} e_j$, for $i\not= j,$ have odd degrees.

\item The degree $0$ part of $A^{n-k,k}$ is the product of rings
$\mathbb{Z}$, one for each element of $B^{n-k,k}$
$$(A^{n-k,k})^0\cong \prod_{a\in B^{n-k,k}} \mathbb{Z} 1_a$$

\begin{figure}[ht!]
\begin{center}
\psfrag{1}{\small$i_1$} \psfrag{2}{\small$i_2$}
\psfrag{3}{\small$i_3$} \psfrag{4}{\small$i_4$}
\epsfig{figure=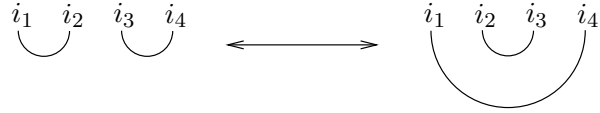} \caption{The elementary change
in $B^{n-k,k}$.} \label{ElementaryMove.figure}
\end{center}
\end{figure}

\end{itemize}

\section{Flat tangles and bimodules}

Denote by $\widehat{B}^m_n$ the space of flat tangles with $m$ top
endpoints and $n$ bottom endpoints (note that in
\cite{Kh2} $\widehat{B}^m_n$ denotes
the space of flat tangles with $2m$ top endpoints and $2n$ bottom
endpoints). Recall that a flat $(m,n)$-tangle $T$ is a proper,
smooth embedding of $\frac{n+m}{2}$ arcs and a finite number of
circles into $\mathbb{R} \times [0,1]$ such that:
\begin{itemize}
\item The boundary points of arcs map to
\begin{equation*}
\{1, 2, . . ., n\}\times \{0\}, \{1, 2, . . ., m\} \times\{1\}.
\end{equation*}
\item Near the endpoints, the arcs are perpendicular to the
boundary of $\mathbb{R}\times [0,1]$.
\end{itemize}
We impose these conditions to ensure that the concatenation of two
such embeddings is still a smooth embedding. Flat tangles
constitute a category with objects -- nonnegative integers, and
morphism from $n$ to $m$ being the isotopy classes of flat $(m,n)$-tangles.\\

\begin{figure}[ht!]
\begin{center}
\psfrag{a}{\small $c$} \psfrag{t}{\small $T$} \psfrag{b}{\small
$b$}\psfrag{atb}{\small $W(c)Tb$} \psfrag{c}{\small closure of
$W(c)Tb$} \psfrag{v}{\small Vertical lines added}
 \epsfig{figure=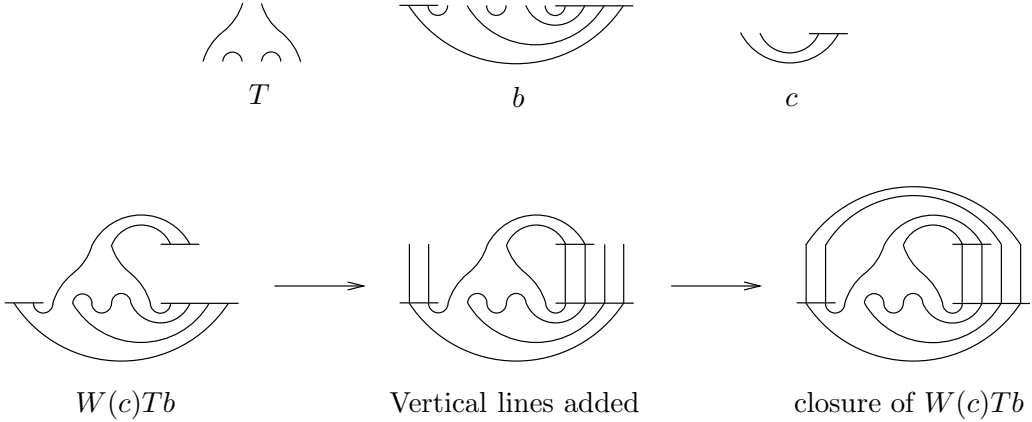}  \caption{Closing $W(c)Tb$.} \label{close.figure}
\end{center}
\end{figure}

\begin{figure}[ht!]
\begin{center}
\psfrag{a}{$c$} \psfrag{T}{$T$}
\psfrag{b1}{$b_1$}\psfrag{F}{$\cF$} \psfrag{b2}{$b_2$}
\psfrag{b3}{$b_3$}
\psfrag{v1}{$0$}\psfrag{v2}{$\mathbb{Z}\{1\}$}\psfrag{v3}{$\cA\{1\}$}
 \epsfig{figure=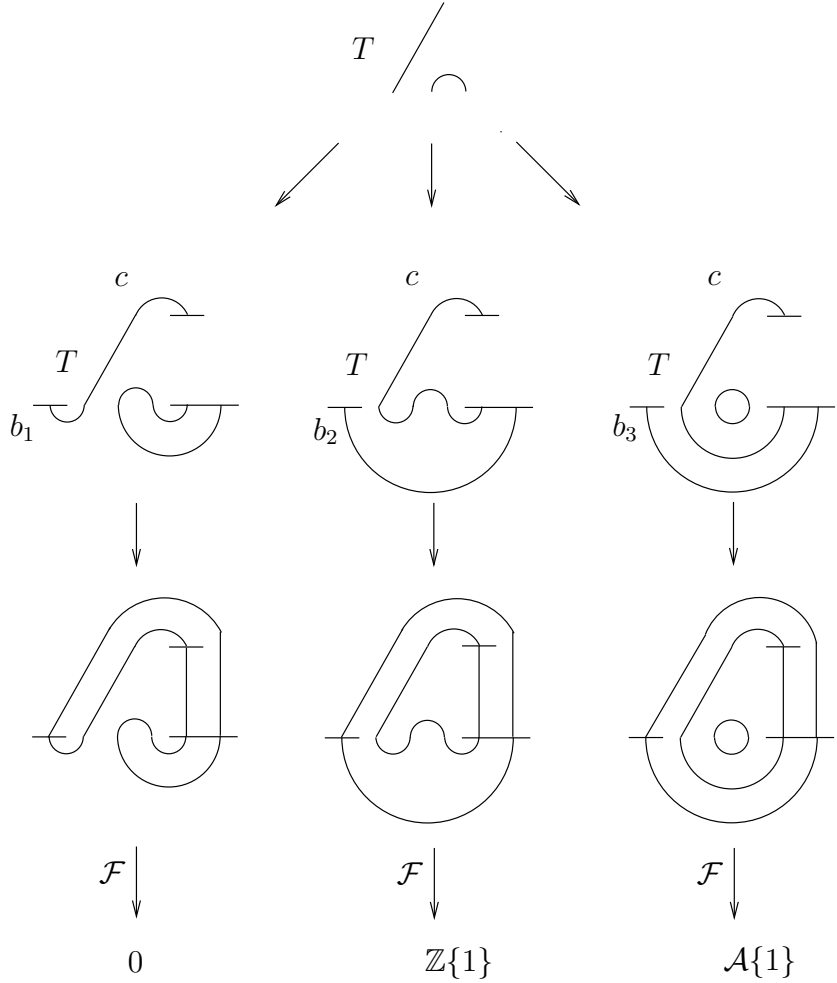}  \caption{Closures of a flat tangle $T$ with fixed sizes
of platforms.} \label{fa.figure}
\end{center}
\end{figure}

To a tangle $T\in \widehat{B}^m_n$ we would like to assign an
$(A^{m-k-l,k+l},A^{n-k,k}$)-bimodule, for all $k$ in the range
$\mathrm{max}(0,\frac{n-m}{2})\leq k \leq
\mathrm{min}(n,\frac{m+n}{2})$ and $l=\frac{m-n}{2}$. First,
define a graded
$(\widetilde{A}^{m-k-l,k+l},\widetilde{A}^{n-k,k})$-bimodule
$\widetilde{\cF}(T)$ by
  \begin{equation*}
    \widetilde{\cF}(T) = \oplusop{b,c} \hspace{0.05in}{_c\widetilde{\cF}(T)_b},
  \end{equation*}
where $b$ ranges over elements of $B^{n-k,k},$ $c$ over elements
of $B^{m-k-l,k+l},$ and
  \begin{equation} \label{def-bimod}
    {_c}\widetilde{\cF}(T)_b \stackrel{\mbox{\scriptsize{def}}}{=} \cF( W(c) T b) \{ n
    \}.
  \end{equation}
Note that $W(c) T b$ is not a union of circles. We close it in the
obvious way and still call it $W(c) T b$ (see
figure~\ref{close.figure}), then apply the functor $\cF$. See
figure~\ref{fa.figure} for an example where $T\in
\widehat{B}^1_3$, $k=2$ and $l=-1,$  $b_i\in B^{1,2}$ and $c\in
B^{0,1}$).
The left action $\widetilde{A}^{m-k-l,k+l} \times
\widetilde{\cF}(T)  \rightarrow \widetilde{\cF}(T)$ comes from
maps
  \begin{equation*}
    \cF(W(a)c) \otimes {_c\widetilde{\cF}(T)_b} \rightarrow
    {_a\widetilde{\cF}(T)_b}.
  \end{equation*}
Likewise, the right action $\widetilde{\cF}(T) \times
\widetilde{A}^{n-k,k} \rightarrow \widetilde{\cF}(T)$ comes from
maps
  \begin{equation*}
    {_c\widetilde{\cF}(T)_b} \otimes \cF(W(b)a)\rightarrow
    {_c\widetilde{\cF}(T)_a}.
  \end{equation*}
Both maps are induced by the obvious ``minimal cobordism'' (see
figure~\ref{contraction.figure}).

Similar to the definition of ${_b(I^{n-k,k})_a}$, we define a
subgroup ${_bI(T)_a}$ of ${_b\widetilde{\cF}(T)_a}$ as follows: If
$W(b)Ta$ contains a type III arc, set
${_bI(T)_a}={_b\widetilde{\cF}(T)_a}$. Otherwise, assuming that
$\cF(W(b)Ta)\cong \cA^{\otimes r}$ in which type II circles
correspond to the first $i$ tensor factors, ${_bI(T)_a}$ is
spanned by elements
 \begin{equation*}
u_1 \otimes \cdots \otimes a_{j-1} \otimes X \otimes u_{j+1}
\otimes \cdots \otimes u_{r} \in \cF(W(b)Ta) \cong \cA^{\otimes
r},
 \end{equation*}
  where $1\leq j\leq i$ and $u_s\in \{\mo, X\}$ for each $1\leq s\leq r$, $s\neq j$.  By taking the direct sum we get a subgroup
 \begin{equation*}
   I(T) \define \oplusop{a \in B^{m-k-l,k+l}, b\in B^{n-k,k}}\hspace{0.05in}
   {_aI(T)_b}.
 \end{equation*}

\begin{lemma} \label{lemma-subbimodule}
$I(T)$ is a subbimodule of $\widetilde{\cF}(T)$.
\end{lemma}

\begin{figure}[ht!]
\begin{center}
\psfrag{x1}{\tiny $X$} \psfrag{x2}{\tiny $\mo$} \psfrag{cx}{\tiny
$C$} \psfrag{a}{\small $a$} \psfrag{b}{\small $b$}
\psfrag{wa}{\small $W(a)$} \psfrag{T}{\small $T$}
\psfrag{c}{\small $c$}\psfrag{wc}{\small
$W(c)$}\psfrag{wca}{\small $\cF(W(c)a)$}\psfrag{watb}{\small
${_a{I(T)}_b}$}\psfrag{wctb}{\small $_c{I(T)}_b$}
\psfrag{mc}{\footnotesize$M=X$} \psfrag{m}{$m$}
\psfrag{delta}{$\Delta$}
\psfrag{in}{\begin{rotate}{-90}$\in$\end{rotate}}\psfrag{ni}{\begin{rotate}{90}$\in$\end{rotate}}
 \epsfig{figure=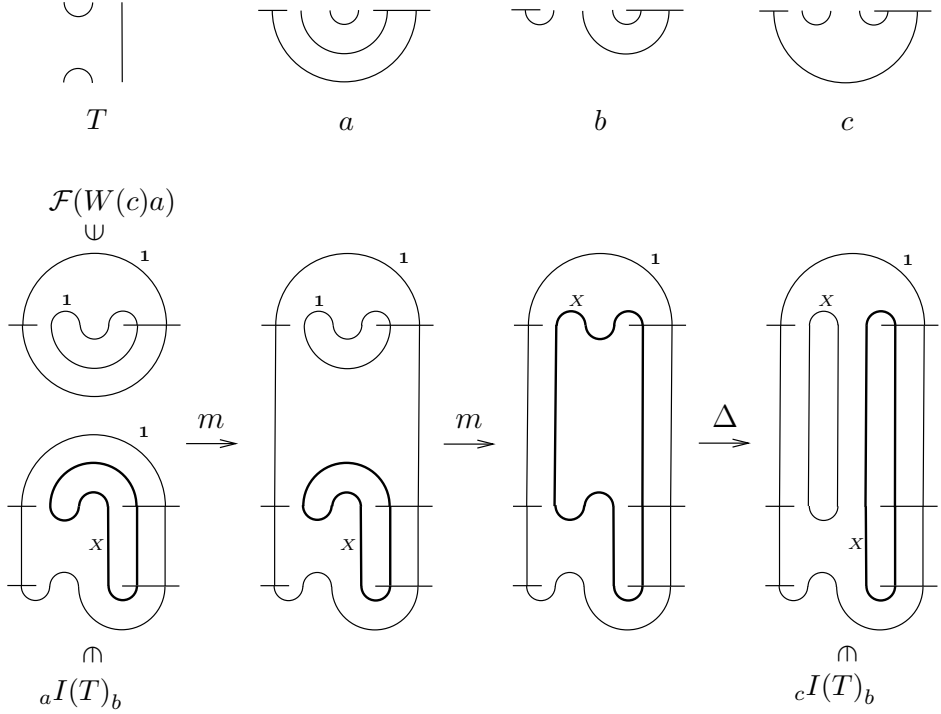}  \caption{Invariance of $I(T)$ under left action.} \label{ideal.figure}
\end{center}
\end{figure}

The proof is similar to that of lemma~\ref{lemma-ideal} and we
omit it. See figure~\ref{ideal.figure} for an example. The
distinguished circle $C$ is thickened.\vspace{0.1in}

Define $\cF(T)$ to be the quotient bimodule of
$\widetilde{\cF}(T)$ over $I(T)$
 \begin{equation*}
   \cF(T) \define \widetilde{\cF}(T) / I(T).
 \end{equation*}

\begin{lemma} \label{lemma-ActionOfIdeal}
The action of $I^{n-k,k}$ on $\cF(T)$ is trivial.
\end{lemma}

The proof is similar to that of lemma~\ref{lemma-ideal}.
\vspace{0.1in}

It follows from the previous lemma that the
$(\widetilde{A}^{m-k-l,k+l},\widetilde{A}^{n-k,k})$-bimodule
structure on $\cF(T)$ descends to an
$(A^{m-k-l,k+l},A^{n-k,k})$-bimodule structure.\\

By taking the direct product over all $0 \leq k \leq n$, we
collect the rings $A^{n-k,k}$ together into a graded ring $A^{n}$
$$A^n\define \prod_{0\leq k\leq n}A^{n-k,k}.$$
As a graded abelian group, $A^n$ is the direct sum of $A^{n-k,k}$,
over $0\leq k\leq n$. Similarily, for a flat tangle $T$, by taking
the direct sum over all $\mathrm{max}(0,\frac{n-m}{2})\leq k \leq
\mathrm{min}(n,\frac{n+m}{2})$ we collect the
$(A^{m-k-l,k+l},A^{n-k,k})$-bimodules $\cF(T)$ into an
$(A^m,A^n)$-bimodule (still call it $\cF(T)$)
$$\cF(T)\define \bigoplus_{\mathrm{max}(0,\frac{n-m}{2})\leq k \leq \mathrm{min}(n,\frac{n+m}{2})}\cF(T).$$
Note that we use the same notation $\cF(T)$ for both
$(A^m,A^n)$-bimodule and individual
$(A^{m-k-l,k+l},A^{n-k,k})$-bimodules.

\begin{prop} Let $T_1,T_2\in \widehat{B}^{m}_{n}$ and $S$ a cobordism between $T_1$ and $T_2$.
Then $S$ induces a degree $\frac{n+m}{2}-\chi(S)$ homomorphism of
$(A^{m},A^{n})$-bimodules
   \begin{equation*}
     \cF(S): \cF(T_1) \to \cF(T_2),
   \end{equation*}
where $\chi(S)$ is the Euler characteristic of $S$.
\end{prop}

\emph{Proof:} We only need to prove the proposition for each
$\mathrm{max}(0,\frac{n-m}{2})\leq k \leq
\mathrm{min}(n,\frac{n+m}{2})$. We have $\widetilde{\cF}(T_1)=
\oplusop{a,b} \cF(W(b)T_1 a) \{ n\}$ and $\widetilde{\cF}(T_2)=
\oplusop{a,b} \cF(W(b)T_2 a) \{ n\}$ where the sum is over $a\in
B^{n-k,k}$ and $b\in B^{m-k-l,k+l}$. The surface $S$ induces a
cobordism $S'=Id_{W(b)}\ S\ Id_a$ from $W(b)T_1 a$ to $W(b)T_2 a$
defined as the ``vertical'' composition of the identity cobordism
from $a$ to $a$, cobordism $S$ from $T_1$ to $T_2$, and the
identity cobordism from $W(b)$ to $W(b)$. $S'$ induces a
homogeneous map of graded abelian groups $\cF(W(b)T_1 a) \to
\cF(W(b)T_2 a)$. Summing over all $a$ and $b$ we get a
homomorphism of
$(\widetilde{A}^{m-k-l,k+l},\widetilde{A}^{n-k,k})$-bimodules
   \begin{equation*}
     \widetilde{\cF}(S): \widetilde{\cF}(T_1) \to
     \widetilde{\cF}(T_2).
   \end{equation*}
Split $S$ into the composition of elementary cobordisms
   \begin{equation*}
     S = S_1\circ S_2 \circ \cdots \circ S_j.
   \end{equation*}
The effect of each elementary cobordism is just an application of
$\iota$, $\varepsilon$, $m$ or $\Delta$. We only need to show that
$\widetilde{\cF}(S)$ takes $I(T_1)$ into $I(T_2)$. This follows
from a argument similar to that in lemma~\ref{lemma-ideal}. The
grading assertion follows from (\ref{DegreeOfCobordism}). Finally,
$\widetilde{\cF}(S)$ is independent of the presentation of $S$ as
the product of elementary cobordisms since $\cF$ is a functor.
$\square$ \vspace{0.1in}

\begin{prop} Isotopic (rel boundary) surfaces induce equal bimodule maps.
\end{prop}

\begin{prop} Let $T_1,T_2,T_3\in \widehat{B}^{m}_{n}$ and $S_1$,
$S_2$ be cobordisms from $T_1$ to $T_2$ and from $T_2$ to $T_3$
respectively. Then $\cF(S_2)\cF(S_1)=\cF(S_2\circ S_1)$.
\end{prop}


\begin{prop} For $T_1\in \widehat{B}^{s}_{n}$, $T_2\in
\widehat{B}^{m}_{s}$ there is a canonical isomorphism of
$(A^{m},A^{n})$-bimodules
   \begin{equation*}
     \cF(T_2 T_1)\cong \cF(T_2) \otimes_{A^{s}}\cF(T_1).
   \end{equation*}
\end{prop}

The proofs of the above propositions are similar to those in
\cite{Kh2}, section $2.7$.

\section{Tangles, complexes of bimodules and tangle cobordisms}

First we recall the definition of tangles. An unoriented
$(m,n)$-tangle $L$ is a proper, smooth embedding of
$\frac{n+m}{2}$ arcs and a finite number of circles into
$\mathbb{R}^2 \times [0,1]$ such that:
\begin{itemize}
\item The boundary points of arcs map to
\begin{equation*}
\{1, 2, . . ., n\}\times \{0\} \times\{0\}, \{1, 2, . . .,
m\}\times \{0\} \times\{1\}.
\end{equation*}
\item Near the endpoints, the arcs are perpendicular to boundary
planes.
\end{itemize}
An oriented $(m,n)$-tangle comes with an orientation of each
connected component. Unoriented tangles constitute a category with
objects -- nonnegative integers and morphisms -- isotopy classes
of $(m,n)$-tangles. The composition of morphism is defined as the
concatenation of tangles. Likewise, oriented tangles form a
category with objects--finite sequences of plus and minus signs,
indicating orientations of the tangle near the endpoints.

A plane diagram of a tangle is a generic projection of the tangle
onto $\mathbb{R}\times [0,1]$. Two diagrams are called isotopic if
they can be transformed into each other through generic
projections. Two plane diagrams represent isotopic tangles if and
only if they can be connected by a chain of diagram isotopies and
Reidemeister moves $R1$, $R2$, and $R3$.

\begin{figure}[ht!]
\begin{center}
 \epsfig{figure=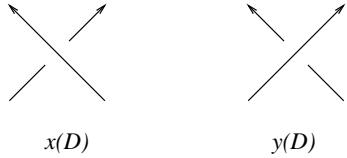}  \caption{Negative and positive crossings.} \label{orient.figure}
\end{center}
\end{figure}

To each diagram $D$ we associate integers $x(D)$ and $y(D)$ which
count the numbers of negative and positive crossings of $D$
respectively, see figure~\ref{orient.figure}.

Fix a diagram $D$ with $s$ crossings of an  oriented
$(m,n)$-tangle $L$. We inductively define the complex of $(A^m,
A^n)$-bimodules $\cF(D)$ associated to $D$ as follows. If $D$ is
crossingless, $\cF(D)$ is the complex with the only nontrivial
term in cohomological degree zero, which is given by the
construction of the previous section.

\begin{figure}[ht!]
\begin{center}
\psfrag{0}{\tiny $0$-smoothing} \psfrag{1}{\tiny $1$-smoothing}
\psfrag{s}{\tiny Saddle cobordism}
 \epsfig{figure=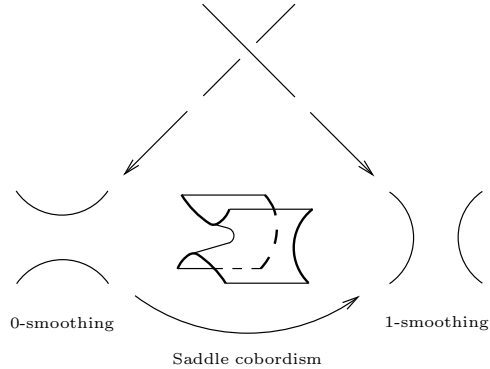}  \caption{Two smoothings of a crossing.} \label{smoothing.figure}
\end{center}
\end{figure}

If the diagram contains one crossing,
consider  the complex $\overline{\cF}(D)$ of $(A^{m},A^{n})$-bimodules
   \begin{equation*}
     0 \to \cF(D(0)) \stackrel{\mbox{\scriptsize{$\partial$}}}{\rightarrow}
\cF(D(1))\{-1\} \to
     0 \label{complex.equation}
   \end{equation*}
where $D(i), i=0,1$ denotes the $i$-smoothing of the crossing
(they are flat $(m,n)$-tangles), $\partial$ is induced by the
obvious ``saddle'' cobordism (see figure~\ref{smoothing.figure}),
and $\cF(D(0))$ sits in the cohomological degree zero.

Inductively, to a diagram with $t+1$ crossings we associate the total
complex $\overline{\cF}(D)$ of the bicomplex
   \begin{equation*}
   0 \to \cF(D(c_0)) \stackrel{\mbox{\scriptsize{$\partial$}}}{\rightarrow}
    \cF(D(c_1))\{-1\} \to 0 \label{complex1.equation}
   \end{equation*}
where $D(c_i), i=0,1$ denotes the $i$-smoothing of a crossing $c$
of $D.$ Finally, define $\cF(D)$ to be $\overline{\cF}(D)$
shifted by $[x(D)]\{ 2x(D)-y(D)\}.$

\begin{theorem} If $D_1$ and $D_2$ are diagrams of an oriented
$(m,n)$-tangle $L$, the complexes $\cF(D_1)$ and $\cF(D_2)$ of
graded $(A^{m},A^{n})$-bimodules are chain homotopy equivalent.
\end{theorem}

Since isotopies of tangles do not involve platforms, the proof of
the theorem is essentially the same as in \cite{Kh2}. It follows
from the above theorem that the isomorphism class of the
complex $\cF(D)$ is an invariant of $L$, denoted by $\cF(L)$. \\

For a graded ring $R$ denote by $\cK(R)$ the category of bounded
complexes of graded $A$-bimodules up to homotopies of complexes.
Objects of $\cK(R)$ are bounded complexes of graded $A$-bimodules
and morphisms of $\cK(R)$ are grading-preserving morphisms of
complexes quotient by null-homotopic ones. We call $M\in \cK(R)$
invertible if there exists $N\in \cK(R)$ such that $N\otimes_R
M\cong R$ and $M\otimes_R N\cong R$ in $\cK(R)$. Here $R$ denotes
the complex $(0\rightarrow R\rightarrow 0)$ with $R$ in
cohomological degree zero. For example, if $L$ is any $n$-stranded
braid, $\cF(L)\in\cK(A^n)$ is invertible. If $M$ is invertible
then
\begin{equation*}
  \mathrm{Hom}_{\cK(R)}(M,M)\cong Z_0(R),
\end{equation*}
where  $Z_0(R)$ is the degree zero component of the center of $R$
(see \cite{Kh2}). Furthermore, we have
\begin{equation*}
  \mathrm{Aut}_{\cK(R)}(M)\cong Z^*_0(R),
\end{equation*}
where $\mathrm{Aut}_{\cK(R)}(M)$ is the group of automorphisms of
$M$ in $\cK(R)$ and $Z^*_0(R)$ is the group of invertible elements
in $Z_0(R)$.\\

For the ring $A^{n-k,k}$ we have

\begin{prop} \label{only-inv}
The only invertible degree $0$ central elements in $A^{n-k,k}$ are
$\pm 1$
\begin{equation*}
  Z^*_0(A^{n-k,k})\cong \{\pm 1\}.
\end{equation*}
\end{prop}
\emph{Proof:} Degree zero elements of $A^{n-k,k}$ have the form
$$v=\sum_{a\in{B^{n-k,k}}} v_a 1_a,$$
where $v_a\in \mathbb{Z}$. For any $a,b\in B^{n-k,k}$ such that $
_a(A^{n-k,k})_b\neq 0$, pick non-zero $x\in \cF(W(b)a)$. Then
$vx=v_a x$ and $xv=v_b x$. If $v$ is central we get $v_a=v_b$. We
can connect any pair $c,d\in B^{n-k,k}$ by a sequence $c=c_0, c_1,
\cdots, c_m=d$ such that $W(c_i)c_{i+1}$ contains no type III
circles. This is equivalent to $_{c_i}(A^{n-k,k})_{c_{i+1}}\neq
0$, so $v_{c_i}=v_{c_{i+1}}$ and $v_c=v_d$. Since $v_a=m$ for all
$a\in B^{n-k,k}$ and some integer $m$, $v=m \sum 1_a=m 1=m$. The
proposition follows.
$\square$\\

From here on we assume familiarity with \cite{Kh3}, where to an
oriented  tangle cobordism there was associated a homomorphism of
complexes of graded $(H^m, H^n)$-bimodules, in a consistent way so
as to produce a projective 2-functor from the 2-category of tangle
cobordisms to the 2-category of natural transformations between
exact functors between homotopy categories of complexes of graded
$H^n$-modules. The construction there extends without difficulty
to our framework. A tangle cobordism can be presented by a movie
$S,$ which is a sequence of Reidemeister moves and critical point
moves. To each consequent pair of tangle diagrams $D_1, D_2$ in a
movie there is associated a natural homomorphism $\cF(D_1)
\longrightarrow \cF(D_2)$ between the corresponding complexes. In
the case of a Reidemeister move, the homomorphism is an
isomorphism in the homotopy category, while for the critical point
moves the homomorphism is induced by either the unit, counit,
multiplication, or comultiplication map on the ring $\mathcal{A}.$

The composition of these homomorphisms gives us a homomorphism
$\cF(S): \cF(D) \longrightarrow \cF(D')$ where $D$ and $D'$ are the first and the
last frame in the movie $S.$ The same argument as in \cite{Kh3}
shows that $\cF(S)= \pm \cF(\tilde{S}),$ where $\tilde{S}$
is any movie between $D$ and $D'$ representing the same cobordism as $S.$
The proposition~\ref{only-inv} above is a necessary ingredient in this argument.

The choice of sign in the equation $\cF(S)= \pm \cF(\tilde{S})$ does not
depend on the sizes of platforms, since the rings $A^{n-k,k}$ are subquotients
of $H^n,$ our bimodules $\cF(D)$ are subquotients of the bimodules in
\cite{Kh2}, and our bimodule homomorphism are induced by those
in \cite{Kh2, Kh3} via subquotient maps. Therefore, the sign is always the
same as in the invariant constructed in \cite{Kh3} and does not depend on
the choice of $k$ between $0$ and $n.$

We can summarize the properties of our construction as follows.

\begin{prop} Complexes $\cF(T)$ of bimodules and homomorphisms
$\pm \cF(S)$ assigned to diagrams of tangle cobordisms assemble into
a projective 2-functor from the 2-category of oriented tangle cobordisms to
the 2-category of natural transformations between exact functors between
homotopy categories of complexes of graded $A^n$-modules.
\end{prop}

Our invariant of tangles and tangle cobordisms carries the same amount of
information as the one in \cite{Kh2, Kh3}. Indeed, since the invariants coming
from the rings $A^n$ are subquotients of the invariants built from $H^n,$
we don't gain new information. On the other hand, the ring $A^{n,n}$ contains
$H^n$ as a subring, since $H^n$ is isomorphic to the direct sum of $\cF(W(b)a)$
over all pairs $a,b$ of diagrams in $B^{n,n}$ which contain $n$
parallel arcs connecting $n$ points on the left platform and $n$ points on
the right platform. The inclusion $H^n \subset A^{n,n}\subset A^{2n}$ extends
to bimodules and bimodule homomorphisms in the two descriptions of
tangle and tangle cobordism invariants. Therefore, the second construction, via
$A^n,$ contains at least as much information as the original one, and our
claim follows.

\section{The Grothendieck group of $A^n$}

\begin{figure}[ht!]
\begin{center}
 \epsfig{figure=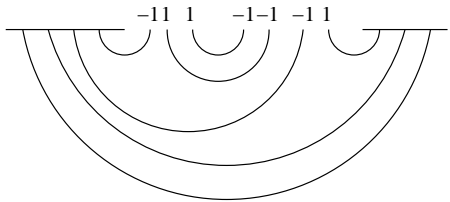}
\caption{Converting $a$ to a sequence of $1$'s and $-1$'s.}
\label{onesm.figure}
\end{center}
\end{figure}

The disjoint union of sets $B^{n-k,k},$ as $k$ ranges from $0$ to $n,$
can be naturally identified with the set $J_n$ of length $n$ sequences
of $1$'s and $-1$'s. A element $a\in B^{n-k,k}$ consists of $n$ arcs
with $2n$ endpoints, $n$ of which lie on platforms and the other $n$
directly between the platforms. We call the endpoints of the second type
\emph{free} endpoints. To each free endpoint we
assign $1$ or $-1$ as follows (see figure~\ref{onesm.figure} for an example).
First, assign $1$ to the left endpoint of each arc and $-1$ to the right
endpoint. We get a sequence of length $2n$ with $n$ ones and
$n$ minus ones. Remove the first $n-k$ and the last $k$ terms in the
sequence (notice that the first $n-k$ terms are all ones, and the last
$k$ are all minus ones). The result is a sequence of length $n$ with
$k$ ones and $n-k$ minus ones. We denote this sequence by $s(a).$

Let $V^n$ be the free $\Z[q,q^{-1}]$-module of rank $2^n$
with the basis $v_s, $ over all sequences $s\in J_n.$ For a sequence
$s=(s_1, \dots, s_n)$ we write $v_s= v_{s_1}\otimes \dots \otimes v_{s_n}$
and identify $V^n$ with the $n$-th tensor power of
the rank $2$ module $V^1.$ Define the weight $w(s)= s_1+ s_2+\dots + s_n$
and $V^n(m)$  to be the subspace of $V_n$ spanned by vectors $v_s$
with $s$ of weight $m.$ Then
$$V^n = \oplusoop{k=0}{n} V^n(2k-n).$$
To each $a\in B^{n-k,k}$ we associated $s(a)\in J_n$ of weight $2k-n.$
We will also denote $v_{s(a)}$ simply by $v_a.$

To each $a\in B^{n-k,k}$ we associate an element $p_a\in V^n$ as follows.
Convert each arc in $a$ disjoint from the platforms into
$$ v_1 \otimes v_{-1} + q v_{-1} \otimes v_1,$$
the indices placed in appropriate positions in the $n$-fold tensor product.
An arc with one end on the left platform and one free end is converted
into $v_{-1},$ in the corresponding position in the tensor product.
An arc with one end on the right platform and one free end
contributes $v_1$ to the tensor product. For example, for $a$
in figure~\ref{onesm.figure},
\begin{eqnarray*}
p_a & = & v_{-1} \otimes v_1 \otimes v_1 \otimes v_{-1}\otimes v_{-1} \otimes
 v_{-1}\otimes v_1 +  \\
   &  &   q \hspace{0.03in}
      v_{-1} \otimes v_1 \otimes v_{-1} \otimes v_{1}\otimes v_{-1} \otimes
      v_{-1}\otimes v_1  + \\
     & &  q\hspace{0.03in}
    v_{-1} \otimes v_{-1} \otimes v_1 \otimes v_{-1}\otimes v_{1} \otimes
   v_{-1}\otimes v_1 +  \\
   & &  q^2 \hspace{0.03in}
      v_{-1} \otimes v_{-1} \otimes v_{-1} \otimes v_{1}\otimes v_{1} \otimes
      v_{-1}\otimes v_1.
\end{eqnarray*}
Notice that
$$p_a = v_a +\hspace{0.1in} \mathrm{ lower}\hspace{0.1in}\mathrm{ order }
\hspace{0.1in} \mathrm{ terms} ,$$ with respect to the order
induced by the relation $1>-1$  ($v_s> v_t$ if $s_i > t_i$ for the
first $i$ where the sequences differ). Hence, $\{p_a\},$ over all
$a\in \sqcup_{k=0}^n B^{n-k,k},$ is a basis of the free
$\Z[q,q^{-1}]$-module $V^n.$

The projective Grothendieck group $K_p(A^n-\mathrm{gmod})$ of the
category of finitely-generated graded projective $A^n$-modules has
generators $[P],$ where $P$ is a projective object of
$A^n-\mathrm{gmod}$ and relations $[P_1]=[P_2]+[P_3]$ whenever
$P_1\cong P_2 \oplus P_3.$ The grading shift functor induces a
$\Z[q,q^{-1}]$-module structure on $K_p(A^n-\mathrm{gmod}).$ An
argument similar to the one in \cite{Kh2} proposition~2 shows that
$P_a\{i\}$ are the only projective indecomposable graded
$A^n$-modules and that $K_p(A^n-\mathrm{gmod})$ is a free
$\Z[q,q^{-1}]$-module of rank $2^n$ with a basis $[P_a],$ $a\in
\sqcup_{k=0}^n B^{n-k,k}.$

Consider the isomorphism of $\Z[q,q^{-1}]$-modules
\begin{equation} \label{eq-iso}
 K_p(A^n-\mathrm{gmod}) \cong V^n
\end{equation}
that takes $[P_a]$ to $p_a.$ For each $(m,n)$-tangle $T$ the
complex of bimodules $\cF(T)$ consists of right projective
bimodules, and the tensor product with $\cF(T)$ is an exact
functor from the category $\mathcal{K}(A^n-\mathrm{gmod})$ to
$\mathcal{K}(A^m-\mathrm{gmod}).$ Here $\mathcal{K}(\mathcal{W})$
denotes the category of bounded complexes of objects of an abelian
category $\mathcal{W}$ up to chain homotopies.

This functor takes a projective object of $A^n-\mathrm{gmod}$ to a complex
of projective objects of $A^m-\mathrm{gmod},$ and hence induces
a homomorphism $[\cF(T)]$ of $\Z[q,q^{-1}]$-modules
$$ K_p(A^n-\mathrm{gmod}) \longrightarrow  K_p(A^m-\mathrm{gmod}).$$
It's easy to compute these maps directly and check that under the
isomorphism (\ref{eq-iso}) they give the standard actions of the
category of tangles on tensor powers
$$V_1^{\otimes n}\cong V^n\otimes_{\Z[q,q^{-1}]} \mathbb{C} $$
of the fundamental representation of quantum $\mathfrak{sl}_2.$

Under this isomorphism the basis of $V^n$ given by
images of indecomposable projective modules $[P_a]$ goes to the Lusztig
 dual canonical basis of $V_1^{\otimes n},$ after changing $q$ to $-q^{-1}$
(the latter basis was explicitly computed in \cite{FK}).

\noindent
Y.~Chen,  Department of Mathematics, University of California, Berkeley, CA 94720.

\noindent
yfchen@math.berkeley.edu

\vspace{0.1in}

\noindent
M.Khovanov, Department of Mathematics, Columbia University, New York, NY 10027.

\noindent
khovanov@math.columbia.edu


\vspace{0.1in}









\begin{thebibliography}{10}

\bibitem{BN}
D.~Bar-Natan,
\newblock Khovanov's homology for tangles and cobordisms,
\emph{Geom. Topol.} 9 (2005) 1443--1499, {\em math.GT/0410495.}

\bibitem{Braden}
T.~Braden,
\newblock Perverse sheaves on Grassmannians,
{\em Canad. J. Math.}, Vol 54, No. 3:493-532, 2002, {\em math.AG/9907152.}

\bibitem{BFK}
J.~Bernstein, I.~Frenkel, and M.~Khovanov,
\newblock A categorification of the Temperley-Lieb algebra and Schur quotients of
U(sl(2)) via projective and Zuckerman functors.
\newblock {\em Selecta Mathematica.}, New ser. 5 : 199--241, 1999.
\newblock {\em math.QA/0002087.}

\bibitem{YC}
Y.~Chen,
\newblock Categorification of level two representations of quantum $sl_N$ via generalized arc rings,
In preparation.

\bibitem{Jones}
V.~Jones,
\newblock
A polynomial invariant for knots via von Neumann
algebras, \emph{Bull. Amer. Math. Soc.} 12 (1985), 103--111.

\bibitem{Kh1}
M.~Khovanov,
\newblock A categorification of the Jones polynomial, \emph{Duke Math. J.}
 101 (2000), no. 3, 359--426, {\em math.QA/9908171.}

\bibitem{Kh2}
M.~Khovanov,
\newblock A functor-valued invariant of tangles,
{\em Algebr. Geom. Topol.}, v.2  (2002), 665-741, {\em math.QA/0103190.}

\bibitem{Kh3}
M.~Khovanov,
\newblock An invariant of tangle cobordisms, {\em Trans. Amer. Math. Soc.} 358
(2006), 315-327, {\em math.QA/0207264.}


\bibitem{FK}
I.~Frenkel and M.~Khovanov,
\newblock Canonical bases in tensor products and graphical calculus for $U_q(sl_2),$
{\em Duke Math. Journal} 87 (1997), No. 3, 409--480.

\bibitem{KR}
M.~Khovanov and L.~Rozansky,
\newblock Matrix factorizations and link homology, \emph{math.QA/0401268.}

\bibitem{KS}
M.~Khovanov and P.~Seidel,
\newblock Quivers, floer cohomology, and braid group actions,
{\em Journal of AMS} 15 (2001), no. 1, 203--271, {\em math.QA/0006056}.

\bibitem{N}
G.~Naot, On the algebraic structure of Bar-Natan's universal geometric complex and
the geometric structure of Khovanov link homology theories, math.GT/0603347.

\bibitem{St1}
C.~Stroppel,
\newblock Categorification of the Temperley-Lieb category, tangles, and cobordisms
via projective functors, {\em Duke. Math. Journal} 126 (2005), no. 3, 547-596.

\bibitem{St2}
C.~Stroppel,
\newblock TQFT with corners and tilting functors in the Kac-Moody case,
{\em math.QA/0605103.}

\bibitem{St3}
C.~Stroppel,
\newblock Perverse sheaves on Grassmannians, Springer fibres and Khovanov homology,
 \emph{math.RT/0608234.}

\end{thebibliography}
\end{document}